\documentclass[a4paper]{article}
\usepackage[affil-it]{authblk}

\usepackage[english]{babel} 
\usepackage{amsmath}
\usepackage{amsfonts}
\usepackage{amssymb}
\usepackage{mathrsfs}
\usepackage{graphicx}
\usepackage{enumerate}
\usepackage{xcolor}
\usepackage{hyperref}
\usepackage[numbers,square]{natbib}
\usepackage{amsthm}
\usepackage{faktor}
\usepackage{dirtytalk}
\usepackage{fancyhdr}

\usepackage[Sonny]{fncychap}

\usepackage{tikz}

\usepackage{mathtools,extarrows}
\usepackage{titlesec} 

\titleformat{\chapter}[display]
  {\bfseries\Large}
  {\filright\MakeUppercase{\chaptertitlename} \Huge\thechapter}
{1ex}
  {\titlerule\vspace{1ex}\filleft}  
  [\vspace{1ex}\titlerule]
\DeclarePairedDelimiter\floor{\lfloor}{\rfloor}

\newenvironment{pf}{Proof:}{$\hspace*{\fill}\Diamond$}
\newtheorem{df}{Definition}[section]
\newtheorem{prop}{Proposition}[section]
\newtheorem{cor}{Corollary}[section]
\newtheorem{ex}{Example}[section]
\newtheorem{rmk}{Remark}[section]

\newtheorem{thm}{Theorem}[section]
\newtheorem{lem}{Lemma}[section]

\begin{document}

\title{Orders occurring as endomorphism ring of a Drinfeld module in some isogeny classes of Drinfeld modules of higher ranks}



\author{Sedric Nkotto Nkung Assong
  \thanks{Electronic address: \texttt{sedric.assong@mathematik.uni-kassel.de}}}
\affil{Institute of Mathematics, University of Kassel, Germany}

\date{Dated: \today}

\maketitle

\begin{abstract}
The question we propose to answer throughout this paper is the following: Given an isogeny class of Drinfeld modules over a finite field, what are the orders of the corresponding endomorphism algebra (which is an isogeny invariant) that occur as endomorphism ring of a Drinfeld module in that isogeny class?\\
It is worth mentioning that this question is different from the ones investigated by the authors Kuhn, Pink in~\cite{kuhn2016finding} and Garai, Papikian in~\cite{garai2019computing}. The former authors rather provided an answer to the question, given a Drinfeld module $\phi$, how does one efficiently compute the endomorphism ring of $\phi$? \\ 
The importance of our question resides in the fact that it might be very helpful to better understand isogeny graphs of Drinfeld modules of higher rank ($r \geq 3$) and may be reopen the debate concerning the application to isogeny-based cryptography.\\
We answer that question for the case whereby the endomorphism algebra is a field by providing a necessary and sufficient condition for a given order to be the endomorphism ring of a Drinfeld module. We apply our result to rank $r=3$ Drinfeld modules and explicitly compute those orders occurring as endomorphism rings of rank 3 Drinfeld modules over a finite field.


\end{abstract}
\newpage
\section{Endomorphism rings of Drinfeld modules of higher ranks}
\textbf{Notations}:\\

$
\begin{array}{rl}
A=\mathbb{F}_q[T]: & \text{ Ring of univariate polynomials in $T$ over a finite field } \mathbb{F}_q=\mathbb{F}_{p^*},~p~ prime.\\
k=\mathbb{F}_q(T): & \text{ Rational function field over } \mathbb{F}_q.\\
L: & \text{ Finite $A$-field}.\\
\mathfrak{p}_v: & \text{ (Generator of the) Kernel of the $\mathbb{F}_q$-algebra homomorphism $\gamma$ defining the $A$-field $L$}.\\
v: & \text{the place of $k$ defined by $\mathfrak{p}_v$}\\
\infty: & \text{The place at infinity of $k$}.\\
m : & \text{The degree of $L$ over $A/\mathfrak{p}_v$ i.e. } m=[L : A/\mathfrak{p}_v]   

\end{array}
$

%
\vspace{1cm}

Nkotto proved in~\cite{nkotto2020explicit}, that a general rank $r$ Weil polynomial defining an isogeny class of rank $r$ Drinfeld modules has the form 
$$M(x)=x^{r_1} + a_1x^{r_1-1} + \cdots + a_{r_1-1}x + \mu\mathfrak{p}_v^{\frac{m}{r_2}}$$ 
where $r_1=[k(\pi) : k]$, $r_2=\sqrt{\dim_{k(\pi)}End\phi \otimes_A k}$ and $r=r_1r_2$.\\
Therefore our restriction on the endomorphism algebra (of the corresponding isogeny class) that must be a field, leads to the restriction to isogeny classes defined by Weil polynomials of the form $$M(x)=x^r + a_1x^{r-1} + \cdots + a_{r-1}x + \mu\mathfrak{p}_v^m$$
We aim in this part to prove the following theorem.\\
\begin{thm}\label{endo_thm}
$A=\mathbb{F}_q[T],~ k=\mathbb{F}_q(T)$ and $\mathfrak{p}_v$ is the (generator of the) kernel of the characteristic morphism $\gamma: A \longrightarrow L$ defining the finite $A$-field $L$.\\
$M(x)= x^r + a_1x^{r-1} + \cdots + a_{r-1}x + \mu\mathfrak{p}_v^m \in A[x]$ is a Weil polynomial, where $m=[L: A/\mathfrak{p}_v\cdot A]$. Let $\mathcal{O}$ be an $A$-order of the function field\\ $k(\pi)=k[x]/M(x)\cdot k[x]$. Let $v_0$ be the unique zero of the Frobenius endomorphism $\pi$ in $k(\pi)$ lying over the place $v$ of $k$.\\
$\mathcal{O}$ is the endomorphism ring of a Drinfeld module in the isogeny class defined by the Weil polynomial $M(x)$ if and only if $\mathcal{O}$ contains $\pi$ and $\mathcal{O}$ is maximal at the place $v_0$. 
\end{thm}
 
\noindent Before proving it, let us recall the notions of Tate modules and Dieudonn\'e modules which are very important for the proof.
\subsection{Tate module of a Drinfeld module}
Let $\psi$ be a Drinfeld module over the $A$-field $L$ with $A$-characteristic $\mathfrak{p}_v$. $v$ denotes the place of $k$ associated to the prime $\mathfrak{p}_v$. Let $\omega$ be a place of $k$ different from $v$ and $\mathfrak{p}_{\omega}$ denotes the corresponding prime. $\psi [\mathfrak{p}_{\omega}^n]$ denotes the group of $\mathfrak{p}_{\omega}^n$-torsion points of $\psi$.
\begin{df}
The Tate module of $\psi$ at $\omega$ is defined by the inverse limit\\
$T_{\omega}\psi := \varprojlim \psi [\mathfrak{p}_{\omega}^n]=Hom_{A_{\omega}}\left( k_{\omega}/A_{\omega}, \psi[\mathfrak{p}_{\omega}^{\infty}]\right) \text{ where } \psi[\mathfrak{p}_{\omega}^{\infty}]=\displaystyle\bigcup_{n\geq 1} \psi[\mathfrak{p}_{\omega}^n]$
\end{df}
\begin{rmk}(Recall)~\\
Let $\phi$ and $\psi$ be two isogenous Drinfeld modules defined over the $A$-field $L$. $Hom_L(\phi, \psi)$ denotes the group of isogenies from $\phi$ to $\psi$. Let $u: \phi \longrightarrow \psi$ be an isogeny. If $y \in \phi[\mathfrak{p}_{\omega}^n]$ then $u(y) \in \psi[\mathfrak{p}_{\omega}^n]$.\\
To $u \in Hom_L(\phi, \psi)\otimes A_{\omega}$, corresponds therefore a canonical morphism of $A_{\omega}$-modules $u^* \in Hom_{A_{\omega}}\left(T_{\omega}\phi, T_{\omega}\psi \right)$. 
\end{rmk}
\begin{thm}\label{Tatetheorem1}[Tate, \cite[see theorem 4.12.12]{DGoss}]~\\
Let $\phi$ and $\psi$ be two isogenous Drinfeld modules over the finite $A$-field $L$ as mentioned in the previous remark. Let $G=Gal(\overline{L}/L)$. The canonical map $$ Hom_{L}(\phi, \psi)\otimes A_{\omega} \xlongrightarrow{\sim}  Hom_{A_{\omega}[G]}\left(T_{\omega}\phi, T_{\omega}\psi \right)$$
is a bijection (as morphism of $A_{\omega}$-modules).
\end{thm}
\begin{cor}\label{Tatecorollary1}~~
\begin{itemize}
\item If $\phi=\psi$ then we have the bijection
$$ End_L\phi \otimes A_{\omega} \xlongrightarrow{\sim} End_{A_{\omega}[G]}T_{\omega}\phi$$
\item We denote $V_{\omega}\phi := T_{\omega}\phi \otimes k_{\omega}$. 
$$ End_L\phi \otimes k_{\omega} \cong End_{k_{\omega}[G]}V_{\omega}\phi$$
as $k_{\omega}$-algebras.
\end{itemize}
\end{cor}
\begin{rmk}\label{charpoly}
Let $\pi$ be the Frobenius endomorphism of the Drinfeld module $\phi$. We denote $M(x)$ the minimal polynomial of $\pi$ over $k$. \\ The characteristic polynomial of the action of $\pi$ on the Tate module $T_{\omega}\phi$ is $M(x)^t$ where $t=\dim_{k(\pi)}End\phi\otimes k$. If $t=1$ as it will be the case in the sequel, then $M(x)$ is the characteristic polynomial of the action of $\pi$ on $T_{\omega}\phi$ .
\end{rmk}

\subsection{Dieudonn\'e module of a Drinfeld module}
We want now to discuss what the so-called Tate's theory says when one works at the place $v$ defined by the $A$-characteristic of the Drinfeld module $\phi$ defined over the finite $A$-field $L$.\\
Let us recall that the Tate's theory at the other places $\omega$, strongly relies on the fact that the polynomial $\phi_{\mathfrak{p}_{\omega}^n}(x)$ is separable. That means $\phi[\mathfrak{p}_{\omega}^n]$ (as group scheme) is \'etale. This is not true anymore at the place $v$. That difficulty is overcome by considering the notion of Dieudonn\'e modules. Before moving forward, let us recall the following theorem known as Dieudonn\'e-Cartier-Oda theorem.
\begin{thm}\label{DieudonneCartier}
Let $m \in \mathbb{N}$ and $L$ be a degree $m$ field extension of $A/\mathfrak{p}_v$. Let $K_v$ be the unique degree $m$ unramified extension of the completion field $k_v$ of $k$ at the place $v$. Let $W$ be the ring of integers of $K_v$. Let $F$ and $V$ be indeterminates such that \\
$FV = VF = \mathfrak{p}_v$\\
$F\lambda = \sigma(\lambda)F$ and $\lambda V = V \sigma(\lambda)~~\forall \lambda \in W$\\
where $\sigma: W \longrightarrow W$ is the unique automorphism induced by the Frobenius $\tau^{\deg\mathfrak{p}_v}$ of $L$.\\
There is an anti-equivalence of categories between the category of finite commutative group scheme over $L$ of finite $A/\mathfrak{p}_v$-rank and the category of left $W[F,V]$-modules of finite $W$-length.  
\end{thm}
\begin{rmk}\label{Dieudonne-rmk}~~
\begin{itemize}
\item Given a finite commutative $L$-group scheme $S$ of finite $A/\mathfrak{p}_v$-rank, we denote $D(S)$ the corresponding left $W[F,V]$-module of finite $W$-length. 
\item $D(S)$ is $W$-free and $rank_{A/\mathfrak{p}_v}S = rank_W D(S)$.
\item $W$ is also known as the ring of Witt vectors over the field $L$ and since $L$ is finite (and therefore perfect), $W$ is a discrete valuation ring and $L$ is its residue field. 
\end{itemize}
\end{rmk}
\begin{df}[Dieudonn\'e module at the place $v$]~~\\
Let $\psi$ be a Drinfeld module over the finite $A$-field $L$ with $m=[L: A/\mathfrak{p}_v]$.\\
The Dieudonn\'e module of $\psi$ is defined by the direct limit
$$T_v\psi := \varinjlim D(\psi[\mathfrak{p}_v^n]) $$
where $D(\psi[\mathfrak{p}_v^n])$ is the left $W[F,V]$-module associated to the $L$-group scheme $\psi[\mathfrak{p}_v^n]$ as mentioned in the previous remark. 
\end{df}
\noindent The corresponding Tate theorem is given below.
\begin{thm}\label{Tatetheorem2}[Serre-Tate, \cite[proposition8.2, corollary 8.3, theorem 8.4]{hartl2005uniformizing}]~\\
The canonical map
$$ Hom_L\left( \phi, \psi \right) \otimes A_v \xlongrightarrow{\sim} Hom_{W[F,V]}\left(T_v\psi, T_v\phi \right)$$
is a bijection (as morphism of $A_v$-modules).
\end{thm}
\begin{rmk}\label{rmkWaterhouse}[see \cite{hartl2005uniformizing}]~
\begin{itemize}
\item If $\phi = \psi$ then we have 
$ End\psi \otimes A_v \xlongrightarrow{\sim} End_{W[F,V]}T_v\psi$
\item We denote $V_v\psi = T_v\psi \otimes K_v$. We have  $ End\psi \otimes k_v \xlongrightarrow{\sim} End_{K_v[F,F^{-1}]}V_v\psi$.
\item $T_v\psi / \mathfrak{p}_v^nT_v\psi$ can be identified to $D(\psi[\mathfrak{p}_v^n])$.
\item The $W[F,V]$-module $D(\psi[\mathfrak{p}_v^n])$ can be decomposed into its \'etale and local parts. $D(\psi[\mathfrak{p}_v^n])=D\left(
\psi[\mathfrak{p}_v^n]\right)_{loc} \oplus D\left(\psi[\mathfrak{p}_v^n]\right)_{\text{\'et}}$.\\
Actually the polynomial $\psi_{\mathfrak{p}_v^n}(x)=x^{nh\deg\mathfrak{p}_v} \cdot g_n(x)$ where $g_n(x)$ is a separable polynomial. \\
$D\left(\psi[\mathfrak{p}_v^n]\right)_{loc}=
D\left(\psi[\mathfrak{p}_v^n]_{loc}\right)$  and $D\left(\psi[\mathfrak{p}_v^n]\right)_{\text{\'et}}
=D\left(\psi[\mathfrak{p}_v^n]_{\text{\'et}}\right)$\\
where $\psi[\mathfrak{p}_v^n]_{loc}=Spec\left(\overline{L}[x]/\langle x^{nh\deg\mathfrak{p}_v} \rangle \right)$ and $\psi[\mathfrak{p}_v^n]_{\text{\'et}}=
Spec\left(\overline{L}[x]/\langle g_n(x) \rangle \right)$.
That means the Dieudonn\'e module can also be decomposed as\\ $T_v\psi= \left(T_v\psi \right)_{loc} \oplus \left(T_v\psi \right)_{\text{\'et}}.$ 
\item The Frobenius $\pi$ of $\psi$ acts on $T_v\psi$ via $\pi = F^m$.
\item $F$ (and therefore $\pi=F^m$) acts on the local part $D\left(\psi[\mathfrak{p}_v]\right)_{loc}$ as a nilpotent element and acts on the \'etale part $D\left(\psi[\mathfrak{p}_v]\right)_{\text{\'et}}$ as an isomorphism. 
\end{itemize}
\end{rmk}
\noindent For more details on this part, one can follow \cite[\S 6, 7 and 8]{hartl2005uniformizing}.\\ 
The following dictionnary can be helpful:
\begin{itemize}
\item $Q=\mathbb{F}_q(C)~\leadsto ~ k=\mathbb{F}_q(T)$
\item $\Gamma\left(C', \mathcal{O}_{C'}\right) ~ \leadsto ~ A=\mathbb{F}_q[T]$
\item $z ~ \leadsto ~ \mathfrak{p}_v$
\item $\mathbb{F}_q[[z]] ~ \leadsto ~ A_v $
\item $\mathcal{O}_{S}[[z]] ~ \leadsto ~ W$
\item $\mathcal{O}_{S}[[z]]\left[\frac{1}{z}\right] ~ \leadsto ~ W[V]$
\item Abelian sheaf $\mathcal{F}  ~ \leadsto ~ $ Drinfeld module $\phi$
\item Dieudonn\'e module $\left(\mathcal{\widehat{F}}, F\right)  ~ \leadsto ~ \text{Dieudonn\'e $W[F,V]$-module } T_v\phi $
\end{itemize}

\subsection{Main theorem}
Before giving the main theorem, let us lay the groundwork with the following lemmas and remarks.
\begin{lem}\label{heightDM}
Let $M(x)=x^r + a_1x^{r-1} + \cdots + a_{r-1}x + \mu\mathfrak{p}_v^m$ be a Weil polynomial as described in the previous chapter. \\
The height $h$ (in the sense of \cite[Definition 4.5.8]{DGoss}) of the isogeny class defined by $M(x)$ is the sub-degree of the polynomial $M(x) \mod \mathfrak{p}_v$. That is\\ $M(x) \equiv x^r + a_1 x^{r-1} + \cdots + a_{r-h}x^h \mod \mathfrak{p}_v$.
\end{lem}
\begin{pf}
Let us first of all recall that the height is an isogeny invariant. That means two isogenous Drinfeld modules share the same height.\\
Let $\psi$ be a Drinfeld module in our isogeny class. We recall that the Dieudonn\'e module $T_v\psi$ of $\psi$ is a $W[F,V]$-module and the Frobenius endomorphism $\pi$ acts on it via $\pi=F^m$ as we mentioned before.\\
$\pi=F^m$ acts $W$-linearly on the Dieudonn\'e module $T_v\psi$ with the same characteristic polynomial (in A[x]) as it does as $A_{\omega}$-linear endomorphism of the Tate module $T_{\omega}\psi$ for any $\omega \neq v$ (see \cite[proof of theorem A1.1.1]{chai2013complex} or replacing Tate modules by Dieudonn\'e modules in the proof of theorem 4 in~\cite[page 167]{mumford1974abelian}).\\
But the characteristic polynomial of the action of $\pi$ on the Tate module $T_{\omega}\psi$ is the minimal polynomial $M(x)$ of $\pi$ over $k$ (since $End\phi \otimes k = k(\pi)$ see remark~\ref{charpoly}).\\
Therefore $M(x)$ is also the characteristic polynomial of the action of the Frobenius endomorphism $\pi=F^m$ on the Dieudonn\'e module $T_v\psi$.\\
One gets from there that $M(x) \mod \mathfrak{p}_v$ is the characteristic polynomial of the action of $\pi$ on $T_v\psi / \mathfrak{p}_v T_v\psi= D\left(\psi[\mathfrak{p}_v]\right)$ (see remark~\ref{rmkWaterhouse}).\\
As mentioned in remark~\ref{rmkWaterhouse}, we also know that $D\left(\psi[\mathfrak{p}_v]\right)$ decomposes (via the corresponding group scheme) into its \'etale and local parts i.e.\\ $D\left(\psi[\mathfrak{p}_v]\right)= D\left(\psi[\mathfrak{p}_v]\right)_{loc} \oplus D\left(\psi[\mathfrak{p}_v]\right)_{\text{\'et}}$.\\
Therefore the characteristic polynomial also splits into
$$M(x) \equiv M_{loc}(x) \cdot M_{\text{\'et}}(x) \mod \mathfrak{p}_v $$
where $M_{loc}(x) \mod \mathfrak{p}_v$ (resp.  $M_{\text{\'et}}(x) \mod \mathfrak{p}_v$) is the characteristic polynomial of the action of $\pi$ on the local part $D\left(\psi[\mathfrak{p}_v]\right)_{loc}$ (resp. on the \'etale part $D\left(\psi[\mathfrak{p}_v]\right)_{\text{\'et}}$). That means,\\
$\deg\left(M_{loc}(x)\mod\mathfrak{p}_v\right) = rank_W D\left(\psi[\mathfrak{p}_v]\right)_{loc}$ and\\
 $\deg\left(M_{\text{\'et}}(x)\mod\mathfrak{p}_v\right) = rank_W D\left(\psi[\mathfrak{p}_v]\right)_{\text{\'et}}$\\
But we have by the definition of the height of $\psi$ (see definition~\ref{def-height})\\
\begin{tabular}{rcl}
$\psi_{\mathfrak{p}_v} $ & = & $\tau^{r\deg\mathfrak{p}_v} + \alpha_1\tau^{r\deg\mathfrak{p}_v-1} + \cdots + \alpha_{(r-h)\deg\mathfrak{p}_v}\tau^{h\deg\mathfrak{p}_v}$\\
 & = & $\left(\tau^{(r-h)\deg\mathfrak{p}_v} + \alpha_1\tau^{(r-h)\deg\mathfrak{p}_v-1} + \cdots + \alpha_{(r-h)\deg\mathfrak{p}_v}\tau^0\right)\tau^{h
\deg\mathfrak{p}_v}$
\end{tabular}~\\
with $\alpha_{(r-h)\deg\mathfrak{p}_v} \neq 0$ That is,\\
$\psi_{\mathfrak{p}_v}(x)$=$\left(x^{q^{(r-h)\deg\mathfrak{p}_v}} + \alpha_1x^{q^{(r-h)\deg\mathfrak{p}_v-1}} + \cdots + \alpha_{(r-h)\deg\mathfrak{p}_v}\right)x^{q^{h
\deg\mathfrak{p}_v}}= g(x) \cdot x^{q^{h
\deg\mathfrak{p}_v}}$
where $g(x)$ is a separable polynomial (since $\alpha_{(r-h)\deg\mathfrak{p}_v} \neq 0$) and \\
$\psi[\mathfrak{p}_v]_{\text{\'et}}=Spec\left(\overline{L}[x]/\langle g(x)\rangle \right)$ and $\psi[\mathfrak{p}_v]_{loc}=Spec\left(\overline{L}[x]/\langle x^{q^{h
\deg\mathfrak{p}_v}} \rangle\right)$\\
where $\overline{L}$ is an algebraic closure of $L$.\\
As we have mentioned in remark~\ref{rmkWaterhouse}, $\pi$ acts on $D\left(\psi[\mathfrak{p}_v]\right)_{loc}$ (resp.  $D\left(\psi[\mathfrak{p}_v]\right)_{\text{\'et}}$) as a nilpotent element (resp. as an isomorphism). That means the characteristic polynomial $M_{loc}(x) \mod \mathfrak{p}_v$ is a power of $x$ and the characteristic polynomial $M_{\text{\'et}}(x) \mod\mathfrak{p}_v$ has only non-zero roots (non-zero eigenvalues). In addition, $\deg \left(M_{\text{\'et}}(x) \mod \mathfrak{p}_v \right) = rank_W D\left(\psi[\mathfrak{p}_v]_{\text{\'et}}\right) = r-h$ (see remark~\ref{Dieudonne-rmk}). Therefore 
$$ M(x)\equiv M_{loc}(x) \cdot M_{\text{\'et}}(x) \equiv x^h\left(x^{r-h} + a_1x^{r-h-1} + \cdots + a_{r-h}\right) \mod \mathfrak{p}_v $$  
and the result follows.

\hspace{13cm}
\end{pf}
\begin{cor}\label{Mloc}
Let $M(x)$ be as in the previous lemma.\\
$M_{loc}(x)$ is the irreducible factor of $M(x)$ in $k_v[x]$ that describes the unique zero of $\pi$ in $k(\pi)$ 
\end{cor}
\begin{pf}
\begin{itemize}
\item First of all $M_{loc}(x)$ is an irreducible factor of $M(x)$ in $k_v[x]$. Indeed,\\
if $M_{loc}(x)=f_i(x)\cdot f_j(x) \in k_v[x]$ is a product of two irreducible factors of $M(x)$ in $k_v[x]$, then since $M_{loc}(x) \equiv x^h \mod \mathfrak{p}_v$, $f_i(x)$ and $f_j(x)$ would have a common zero modulo $\mathfrak{p}_v$. That is not possible since $M(x)$ is a Weil polynomial.
\item If $f_{i_0}(x)$ is the factor of $M(x)$ in $k_v[x]$ describing the zero $\mathfrak{p}_{i_0}$ of $\pi$ in $k(\pi)$, then the constant coefficient $a_{0,i_0}$ of $f_{i_0}(x)$ must be divisible by $\mathfrak{p}_v$. Indeed,\\ 
$v_{i_0}(\pi)> 0$ i.e. $\overline{v}\circ \tau_{i_0}(\pi)> 0$. In other words $\overline{v}(\pi_{i_0}) > 0$,\\ 
where $\pi_{i_0}$ denotes a root of $f_{i_0}(x)$.\\
That means, $\overline{v}_{\vert k_v(\pi_{i_0})}\left(\pi_{i_0}\right) > 0$ i.e. $v_{i_0}(\pi_{i_0}) > 0$.
 \\
As a result $v_{i_0}\left(N_{k_v(\pi_{i_0})/k_v}\left(\pi_{i_0}\right)\right) > 0$ and thus\\ $v\left( N_{k_v(\pi_{i_0})/k_v}\left(\pi_{i_0}\right)\right)>0
\text{ since } N_{k_v(\pi_{i_0})/k_v}\left(\pi_{i_0}\right) \in k_v.$\\
But the constant coefficient of $f_{i_0}(x) , ~~a_{0,i_0}=(-1)^{\deg f_{i_0}(x)}N_{k_v(\pi_{i_0})/k_v}\left(\pi_{i_0}\right)$. 
That means we also have $v(a_{0,i_0})>0$ and the claim follows.
\item Since $M(x)$ is a Weil polynomial, there must be only one such factor $f_{i_0}(x)$ of $M(x)$ in $k_v[x]$. Since $M_{loc}(x)\equiv x^h \mod \mathfrak{p}_v$, the constant coefficient of $M_{loc}(x)$ in $A_v[x]$ is divisible by $\mathfrak{p}_v$.\\
Hence $M_{loc}(x)=f_{i_0}(x)$ is the irreducible factor of $M(x)$ in $k_v[x]$ describing the zero $\mathfrak{p}_{i_0}$ of $\pi$ in $k(\pi)$.
\end{itemize}
\hspace{13.5cm}
\end{pf}
Before moving forward, let us formulate the problem.\\

\underline{\textbf{Formulation of the problem}}:\\

Yu in~\cite{yu1995isogenies} basically showed that for an isogeny class of rank 2 Drinfeld modules, the orders occurring as endomorphism ring of a Drinfeld module are either (in case the endomorphism algebra is not a field) the maximal orders in the quaternion algebra over $k$ ramified at exactly the places $v$ and $\infty$, or those orders $\mathcal{O}$ of $k(\pi)$ containing $\pi$ that are maximal at all the places lying over $v$ i.e. such that $\mathcal{O} \otimes A_v$ is a maximal $A_v$-order of the $k_v$-algebra $k_v(\pi).$\\ 
Now the question is: What about Drinfeld modules of higher rank ($r \geq 3$)? \\
Of course for an order $\mathcal{O}$ of (the endomorphism algebra) $k(\pi)$ to be the endomorphism ring of a Drinfeld module, it is necessary that the Frobenius $\pi \in \mathcal{O}$. But must we have $\mathcal{O}$ maximal at all the places of $k(\pi)$ lying over the place $v$? In other words, must we have $\mathcal{O} \otimes A_v$ maximal $A_v$-order of the $k_v$-algebra $k_v(\pi)$?  The answer is No! and we provide below an example of a rank 3 Drinfeld module whose endomorphism ring is not at all places of $k(\pi)$ lying over the place $v$ maximal.\\
Before the example, let us recall the definition and a fact concerning the notion of conductor of an order.
\begin{df}[Recall]
$A=\mathbb{F}_q[T], ~ k=\mathbb{F}_q(T)$\\
Let $F/k$ be a function field and $\mathcal{O}_{max}$ be the ring of integers of $F$. Let $\mathcal{O}$ be an $A$-order of $F$. The conductor $\mathfrak{c}$ of $\mathcal{O}$ is the maximal ideal of $\mathcal{O}$ which is also an ideal of $\mathcal{O}_{max}$. It is defined by $\mathfrak{c}= \displaystyle\{ x \in F \mid x\mathcal{O}_{max} \subseteq \mathcal{O} \}$. 
\end{df}

\begin{rmk}
As a very well known fact, $disc\left( \mathcal{O} \right)= N_{F/k}\left(\mathfrak{c}\right)disc\left(\mathcal{O}_{
max}\right)$.\\
Where $disc(?)$ denotes the discriminant of a basis of the corresponding free $A$-lattice and $N_{F/k}(?)$ denotes the norm of the ideal in argument.  We recall that if $\mathfrak{P}$ is a prime  of F above the prime $\mathfrak{p}$ of k then $N_{F/k}\left(\mathfrak{P}\right)=\mathfrak{p}^{\mathfrak{f}}$ where $\mathfrak{f}$ denotes the residual degree of  $\mathfrak{P} \mid \mathfrak{p}.$ In addition $N_{F/k}(?)$ is multiplicative i.e. $N_{F/k}\left(\mathfrak{P}_1\mathfrak{P}_2\right)= N_{F/k}\left(\mathfrak{P}_1\right)
N_{F/k}\left(\mathfrak{P}_2\right)$.
\end{rmk}
\begin{ex}~\\ $A=\mathbb{F}_5[T], ~ k=\mathbb{F}_5(T), ~ L=\mathbb{F}_{125}=\mathbb{F}_5(\alpha) \text{ with } \alpha^3+3\alpha+3=0$.\\
$\mathfrak{p}_v=Ker\gamma = \langle T \rangle$. $M(x)=x^3+(T+1)x^2 +(T^2+3T+4)x + 4T^3$.\\
One easily shows that $M(x)$ is a Weil polynomial (see~\cite{nkotto2020explicit})\\
$disc\left(M(x)\right)=T^2(T+4)^2(T^2+4T+2)$. Following the paper~\cite{scheidler2004algorithmic} one computes the following:\\
The discriminant of the cubic function field $k(\pi)$ is\\ $\Delta = disc\left(k(\pi)\right) = (T+4)^2(T^2+4T+2)$. We set $I=\sqrt{\frac{disc\left(M(x)\right)}{\Delta}}=T$.\\
The maximal order of the function field $k(\pi)/k$ is the order generated by $\langle \omega_0, \omega_1, \omega_2 \rangle$, where $\omega_0=1,~ \omega_1= \tilde{\pi}=\pi + 2T+2, ~ \omega_2 = \frac{\alpha_2 + \beta_2\tilde{\pi} + \tilde{\pi}^2}{I}= \frac{\alpha_2 + \beta_2\tilde{\pi} + \tilde{\pi}^2}{T}$\\
With\\
$
\begin{cases}
3\beta_2^2 + c_1 \equiv 0 \mod I\\
\beta_2^3+ c_1\beta_2 + c_2 \equiv 0 \mod I^2\\
\alpha_2 \equiv -2\beta_2^2 \equiv\frac{2c_1}{3} \mod I
\end{cases}
$ \\
Where $c_1$ and $c_2$ denote the coefficients of the so-called standard form of the cubic polynomial $M(x)$. We will come back later on to this. \\
After solving the system, one gets $\beta_2=4$ and $\alpha_2=3$.\\ 
That is, $\omega_2=\frac{3+ 4(\pi + 2T+2) + (\pi + 2T+2)^2}{T}$ \\
We now claim that the conductor $\mathfrak{c}$ of $\mathcal{O}=A[\pi]$ is $\mathfrak{c}=T\cdot \mathcal{O} + (\pi -3T +3)\cdot \mathcal{O}$.\\
Indeed,
$$M(x) \equiv x(x-3T+3)^2 \mod T$$ 
We also have $(\pi-3T+3)(\lambda_0\omega_0 + \lambda_1\omega_1 + \lambda_2\omega_2) \in A[\pi] \text{ for } \lambda_i \in A$. Because $(\pi-3T+3)\omega_2 = (T+1)\pi + 4T^2+4T +3 \in A[\pi]$. That means $\pi -3T+3 \in \mathfrak{c}$.\\
Therefore $T\cdot \mathcal{O} + (\pi -3T +3)\cdot \mathcal{O} \subseteq \mathfrak{c} \varsubsetneq \mathcal{O}.$\\
Let us consider the canonical morphisms 
$$
\begin{array}{rl}
\displaystyle A[\pi]\simeq A[x]/M(x)\cdot A[x] \xlongrightarrow{~\varphi_1~} \frac{(A/T\cdot A)[x]}{M(x)\cdot (A/T\cdot A)[x])} \xlongrightarrow{~\varphi_2~} & \frac{(A/T\cdot A)[x]}{(x-3T+3)\cdot (A/T\cdot A)[x])} \\
& \simeq A/T\cdot A
\end{array} 
 $$
$T\cdot\mathcal{O} + (\pi - 3T+3)\cdot\mathcal{O}$ is a maximal ideal of $\mathcal{O}$ as kernel of the morphism $\varphi_2\circ\varphi_1$ since $A[\pi]/Ker(\varphi_2\circ\varphi_1) \simeq Im(\varphi_2\circ\varphi_1)\simeq A/T\cdot A$ is a field.
Therefore $\mathfrak{c}= T\cdot \mathcal{O} + (\pi -3T +3)\cdot \mathcal{O}$.\\
$M(x)\equiv x(x+3)^2 \mod T$. Since $M(x)$ is a Weil polynomial, the irreducible decomposition of $M(x)$ over the completion field $k_v$ is of the form \\
$M(x)=M_1(x) \cdot M_2(x) \in k_v[x]$. That means $\mathfrak{p}_v=T$ splits into two primes $\mathfrak{p}_1$ and $\mathfrak{p}_2$ in $k(\pi)$.\\
As a matter of fact, any prime ideal $\mathfrak{p}$ of $\mathcal{O}$ containing $T$ is either\\ $T\cdot\mathcal{O} + (\pi -3T+3)\cdot \mathcal{O}$ or $T\cdot\mathcal{O} + \pi\cdot \mathcal{O}.$ Indeed,\\
First of all $T\cdot\mathcal{O} + (\pi -3T+3)\cdot \mathcal{O}$ and $T\cdot\mathcal{O} + \pi\cdot \mathcal{O}$ are maximal ideals of $\mathcal{O}=A[\pi]$ as kernel of the canonical morphisms 
$$
\begin{array}{rl}
\displaystyle A[\pi]\simeq A[x]/M(x)\cdot A[x] \xlongrightarrow{~\varphi_1~} \frac{(A/T\cdot A)[x]}{M(x)\cdot (A/T\cdot A)[x])} \xlongrightarrow{~\varphi_2~} &  \frac{(A/T\cdot A)[x]}{(x-3T+3)\cdot (A/T\cdot A)[x])}\\
 &  \simeq A/T\cdot A 
\end{array}
 $$
and
$$ 
\begin{array}{rl}
\displaystyle A[\pi]\simeq A[x]/M(x)\cdot A[x] \xlongrightarrow{~\varphi_1'~} \frac{(A/T\cdot A)[x]}{M(x)\cdot (A/T\cdot A)[x])} \xlongrightarrow{~\varphi_2'~} & \frac{(A/T\cdot A)[x]}{x\cdot (A/T\cdot A)[x])} \\
&  \simeq A/T\cdot A 
\end{array}
$$
respectively.\\
Since $M(x)\equiv x(x-3T+3)^2 \mod T $ and $M(\pi)=0$, we have\\
$\pi(\pi-3T+3)^2 \in T\cdot A[\pi] \subseteq \mathfrak{p}$. But $\mathfrak{p}$ is a prime ideal of $\mathcal{O}$. That means $\pi \in \mathfrak{p}$ or $\pi-3T+3 \in \mathfrak{p}$. In other words\\
$T\cdot\mathcal{O} + (\pi -3T+3)\cdot \mathcal{O} \subseteq \mathfrak{p}$ or $T\cdot\mathcal{O} + \pi\cdot \mathcal{O} \subseteq \mathfrak{p}$\\ 
From the maximality of these ideals we conclude that\\ 
$\mathfrak{p}= T\cdot\mathcal{O} + (\pi -3T+3)\cdot \mathcal{O}$ or $\mathfrak{p}= T\cdot\mathcal{O} + \pi\cdot \mathcal{O}.$ \\
We assume then WLOG that $\mathfrak{p}_2 \cap \mathcal{O} = T\cdot\mathcal{O} + (\pi -3T+3)\cdot\mathcal{O}=\mathfrak{c}$.\\ 
That is, $\mathfrak{p}_2 \mid \mathfrak{c}$ and $\mathfrak{p}_1 \nmid \mathfrak{c}$.\\
The norm of the conductor is\\ $N_{k(\pi)/k}\left(\mathfrak{c}\right)=T^2$ since $disc\left(M(x)\right)=T^2\cdot disc\left(k(\pi)\right)$.\\
Therefore we have only two possibilities for orders occurring as endomorphism of a Drinfeld module: $A[\pi]$ and the maximal order $\mathcal{O}_{max}$. This is due to the fact that the norm of the conductor of any order $\mathcal{O}$ containing properly $A[\pi]$ (i.e. $A[\pi] \varsubsetneq \mathcal{O} \subseteq \mathcal{O}_{max}$) is a square of a proper divisor of $T^2$ and thus must be a unit. In other words $disc(\mathcal{O})=disc\left(\mathcal{O}_{max}\right)$. i.e. $\mathcal{O}=\mathcal{O}_{max}.$\\
After some computations (using a code we implemented in the computer algebra system SAGE) we found the following:
\begin{itemize}
\item For $\phi_T = -\alpha^2\tau^3 + 2\alpha^2\tau^2 + \alpha^2\tau $ we have:
$$\omega_2= \frac{3 + 4(\tau^3 + 2\phi_T + 2) + (\tau^3 + 2\phi_T + 2)^2}{\phi_T} \in L\{\tau\} \text{ and }\phi_T \cdot \omega_2 = \omega_2 \cdot \phi_T.$$
In other words $\omega_2 \in End\phi$. Therefore $End\phi = \mathcal{O}_{max}$.
\item For $\psi_T = \tau^3 + \tau^2 + \tau$ we have:
 $$\omega_2= \frac{3 + 4(\tau^3 + 2\psi_T + 2) + (\tau^3 + 2\psi_T + 2)^2}{\psi_T} \notin L\{\tau\} \text{ and a fortiori } \omega_2 \notin End\psi.$$
 Since we have only two possibilities for $End\psi$, we can conclude that $End\psi = A[\pi]$.\\
$A[\pi]$ is therefore the endomorphism ring of a Drinfeld module but $A[\pi]$ is not maximal at at least one of the places of $k(\pi)$ lying over the place $v$ because its conductor $\mathfrak{c}$ is not relatively prime to $\mathfrak{p}_v=T$.  
\end{itemize} 
\end{ex}
\noindent One can notice in the example above that $M_{loc}(x)=M_1(x) \equiv x \mod \mathfrak{p}_v$. That means $\deg M_{loc}(x)=1$. Thus any order containing $\pi$ is maximal at the corresponding place $v_1$ (which represent the zero of $\pi$ in $k(\pi)$).\\
Concerning the \'etale part, \\
$M_{\text{\'et}}(x) = M_2(x) \equiv (x+3)^2 \mod \mathfrak{p}_v$. i.e. $\deg M_{\text{\'et}}(x)=2$.\\
We have then here ``enough" $\mathfrak{p}_v$-torsion points. \\
This example already encodes some tips for the generalization.
\begin{df}\cite[remark 4.7.12.1]{DGoss}[recall]~~\\
Let $\phi$ and $\psi$ be two isogenous Drinfeld modules over $L$. Let $u: \phi \longrightarrow \psi,~ u \in L\{\tau \}$ be an isogeny from $\phi$ to $\psi$.\\
$\psi$ is called the quotient of the Drinfeld module $\phi$ by the kernel $G$ of $u$ and denoted $\psi := \phi/G$.

\end{df}
\begin{lem}\label{lemma1}
Let $\phi$ be a Drinfeld module over the finite $A$-field $L$ whose endomorphism algebra is a field i.e. $End\phi \otimes k= k(\pi)$, where $\pi$ is the Frobenius endomorphism of $\phi$. Let $\mathcal{O}$ be an $A$-order of $k(\pi)$ containing $\pi$. We choose a place $\omega$ of $k$ different from $v$.\\
If $End\phi \otimes A_{\omega} \not\cong \mathcal{O} \otimes A_{\omega}$ as $A_{\omega}$-module then there exists a Drinfeld module quotient $\psi = \phi/G_{\mathcal{L}}$ such that\\
$End\psi \otimes A_{\omega} \cong \mathcal{O} \otimes A_{\omega}$ and $End\psi\otimes A_{\nu} \cong End\phi \otimes A_{\nu}$ for all places $\nu \neq \omega$. 
\end{lem}
\begin{pf}
With the hypotheses of the lemma, \\
let us assume that $End\phi \otimes A_{\omega} \not\cong \mathcal{O} \otimes A_{\omega}$. We are looking for an isogeny $u$ that changes (via its kernel) the Drinfeld module $\phi$ into a Drinfeld module $\psi$ so that the endomorphism ring of the resulting Drinfeld module coincides at $\omega$ with $\mathcal{O}$.\\
$\mathcal{O}$ is an $A$-order of $k(\pi)$ containing $\pi$. That means $\mathcal{O} \otimes A_{\omega}$ is an $A_{\omega}$-order of the $k_{\omega}$-algebra $k_{\omega}(\pi)=End\phi \otimes k_{\omega}$. We also know from the corollary~\ref{Tatecorollary1} of the Tate theorem that there is a canonical isomorphism of $k_{\omega}$-algebras $End\phi \otimes k_{\omega} \xlongrightarrow{\sim} End_{k_{\omega}[\pi]}V_{\omega}\phi$, where $V_{\omega}\phi = T_{\omega}\phi \otimes k_{\omega}$.\\
Since in addition $\pi \in \mathcal{O}$, $V_{\omega}\phi$ therefore contains an $A_{\omega}$-lattice $\mathcal{L}$ containing $T_{\omega}\phi$ and stable under the action of $\pi$ such that the corresponding order $End_{A_{\omega}[\pi]}\mathcal{L} \cong \mathcal{O}\otimes A_{\omega}$ as $A_{\omega}$-modules. We consider then such an $A_{\omega}$-lattice $\mathcal{L}$. We have then $T_{\omega}\phi \subseteq \mathcal{L} \subseteq V_{\omega}\phi$.\\
Let $(t_1; \cdots, t_r)$ be an $A_{\omega}$-basis of $T_{\omega}\phi$ and $(z_1, \cdots, z_r)$ be an $A_{\omega}$-basis of $\mathcal{L}$, where $r=rank \phi$. $M_0$ denotes the matrix in $\mathscr{M}_{r \times r}\left( A_{\omega} \right)$ such that
$$
\begin{pmatrix}
t_1 \\ \vdots \\ t_r
\end{pmatrix}=M_0
\begin{pmatrix}
z_1\\ \vdots \\ z_r
\end{pmatrix}
$$
Let $s=\omega\left(det M_0 \right)$ be the valuation (wrt $\omega$) of the determinant $det M_0$.\\
$det M_0= \alpha_0\mathfrak{p}_{\omega}^s$, where $\mathfrak{p}_{\omega}$ is the uniformizing element of the place $\omega$ and $\alpha_0$ is a unit in $A_{\omega}$. The reader can notice that $s>0$ because $End\phi \otimes A_{\omega} \not\cong \mathcal{O} \otimes A_{\omega}$.\\
We consider the following map
$$
\begin{tabular}{rcc}
$Co(M_0)^t: T_{\omega}\phi$ & $\xlongrightarrow{\hspace{1cm}}$ & $\mathcal{L}$\\
$
\begin{pmatrix}
t_1 \\ \vdots \\ t_r
\end{pmatrix}$ & $\longmapsto$ & 
$\alpha_0\mathfrak{p}_{\omega}^s
\begin{pmatrix}
z_1\\ \vdots \\ z_r
\end{pmatrix}$
\end{tabular}
$$
The kernel of this map is $ker Co(M_0)^t=M_0\cdot \phi[\mathfrak{p}_{\omega}^s]$.\\
We recall that $Co(M_0)^t$ (as one can guess) denotes the transpose of the co-matrix of the matrix $M_0$.\\ 
Indeed, \\
if $\lambda_1t_1 + \cdots + \lambda_rt_r \in M_0\cdot \phi[\mathfrak{p}_{\omega}^s]$ then\\ 
$Co(M_0)^t\cdot \left( \lambda_1t_1 + \cdots + \lambda_rt_r \right) \in Co(M_0)^t\cdot M_0\cdot \phi[\mathfrak{p}_{\omega}^s]=\mathfrak{p}_{\omega}^s \cdot \phi[\mathfrak{p}_{\omega}^s]=\{0\}$.\\
That is, $Co(M_0)^t\cdot \left( \lambda_1t_1 + \cdots + \lambda_rt_r \right)=0$ and thus\\ $\lambda_1t_1 + \cdots + \lambda_rt_r \in Ker Co(M_0)^t$. \\
Conversely if $\lambda_1t_1 + \cdots + \lambda_rt_r \in Ker Co(M_0)^t$ then $Co(M_0)^t\cdot \left( \lambda_1t_1 + \cdots + \lambda_rt_r \right)$=$0$ \\
i.e. $\alpha_0\mathfrak{p}_{\omega}^s\left( \lambda_1z_1 + \cdots + \lambda_rz_r \right)=0$ and therefore $\lambda_1z_1 + \cdots + \lambda_rz_r \in \phi[\mathfrak{p}_{\omega}^s]$.\\
That means $\lambda_1t_1 + \cdots + \lambda_rt_r=M_0 \cdot \left(\lambda_1z_1 + \cdots + \lambda_rz_r \right) \in M_0\cdot \phi[\mathfrak{p}_{\omega}^s]$.  \\
Hence  $ker Co(M_0)^t=M_0\cdot \phi[\mathfrak{p}_{\omega}^s]$. \\
Applying the first isomorphism theorem to the morphism of $A_{\omega}$-modules, one gets $T_{\omega}\phi / M_0\cdot \phi[\mathfrak{p}_{\omega}^s] \cong Im\left(Co(M_0)^t\right)=\langle \mathfrak{p}_{\omega}^sz_1, \cdots , \mathfrak{p}_{\omega}^sz_r \rangle$.\\
Let $\mathcal{L}_s=\langle \mathfrak{p}_{\omega}^sz_1, \cdots , \mathfrak{p}_{\omega}^sz_r \rangle$ be the $A_{\omega}$-lattice generated by $\left( \mathfrak{p}_{\omega}^sz_1, \cdots , \mathfrak{p}_{\omega}^sz_r \right)$.\\
$T_{\omega}\phi / M_0\cdot \phi[\mathfrak{p}_{\omega}^s] \cong \mathcal{L}_s = \mathfrak{p}_{\omega}^s\cdot \mathcal{L}$.\\
We set $G_{\mathcal{L}}=M_0\cdot \phi[\mathfrak{p}_{\omega}^s]$ and we consider the Drinfeld module quotient $\psi = \phi / G_{\mathcal{L}}$ defined over $L$.\\
The existence of the Drinfeld module $\psi$ is guaranteed by the fact that the separable additive polynomial 
$$ u = x\displaystyle\prod_{\alpha \in G_{\mathcal{L}}}\left( 1 - \frac{x}{\alpha} \right) $$
whose kernel $G_{\mathcal{L}}$ ( which is stable under the action of the Frobenius endomorphism $\pi$ mainly because $\pi \in \mathcal{O}$), lie in $L\{ \tau\}$ (see \cite[proposition 1.1.5 and corollary 1.2.2]{DGoss}), in addition to the fact that the local part of the group scheme $H=Spec\left(\overline{L}[x]/\langle u(x) \rangle \right)$ is trivial because $u \in L\{\tau\}$ is separable (see \cite[proposition 4.7.11, for t=0]{DGoss}).\\
We have then $T_{\omega}\psi \cong T_{\omega}\phi / M_0\cdot \phi[\mathfrak{p}_{\omega}^s] \cong \mathcal{L}_s = \mathfrak{p}_{\omega}^s\cdot \mathcal{L}$ as $A_{\omega}$-modules.\\
Since $G_{\mathcal{L}}= M_0\cdot \phi[\mathfrak{p}_{\omega}^s]$ and $\mathcal{L}$ are stable under the action of $\pi$, so are $T_{\omega}\psi$ and $\mathcal{L}_s$. In other words $T_{\omega}\psi \cong \mathcal{L}_s$ as $A_{\omega}[\pi]$-modules.\\
That means $End_{A_{\omega}[\pi]}T_{\omega}\psi \cong End_{A_{\omega}[\pi]}\mathcal{L}_s$.\\
One also easily checks that (since $\mathcal{L}_s=\mathfrak{p}_{\omega}^s\cdot \mathcal{L}$) $\mathcal{L}$ and $\mathcal{L}_s$ generate the same order i.e. $End_{A_{\omega}}\mathcal{L}_s = End_{A_{\omega}}\mathcal{L}$.\\
Therefore $End_{A_{\omega}[\pi]}T_{\omega}\psi \cong End_{A_{\omega}[\pi]}\mathcal{L}$. Applying the Tate theorem~\ref{Tatetheorem1}, one gets then $End\psi \otimes A_{\omega} \cong End_{A_{\omega}[\pi]}T_{\omega}\psi \cong End_{A_{\omega}[\pi]}\mathcal{L} \cong \mathcal{O} \otimes A_{\omega}$.\\
At all the other places $\nu \neq \omega, v$ of $k$, we have the following:\\
$0 \longrightarrow G_{\mathcal{L}}=M_0\cdot\phi[\mathfrak{p}_{\omega}^s] \xhookrightarrow{\hspace{1cm}} \phi \xlongrightarrow{\hspace{0.4cm} u\hspace{0.4cm}} \psi \longrightarrow 0$ is an exact sequence.\\
$G_{\mathcal{L}}$ has no non-trivial $\mathfrak{p}_{\nu}$-torsion points. Applying the Tate theorem at the place $\nu$ to this short exact sequence, one gets the exact sequence\\ $0 \longrightarrow T_{\nu}\phi \xlongrightarrow{\hspace{1cm}} T_{\nu}\psi \longrightarrow 0$. That means $T_{\nu}\phi \cong T_{\nu} \psi$ as $A_{\nu}$-modules.\\
In other words $End\psi \otimes A_{\nu} \cong End_{A_{\nu}[\pi]}T_{\nu}\psi \cong  End_{A_{\nu}[\pi]}T_{\nu}\phi \cong End\phi \otimes A_{\nu}$.

\hspace{13cm}
\end{pf}
\begin{lem}\label{lemma2}
Let $\phi$ be a Drinfeld module over the finite $A$-field $L$ whose endomorphism algebra $End\phi \otimes k= k(\pi)$ is a field, where $\pi$ denotes the Frobenius endomorphism of $\phi$. Let $\mathcal{O}$ be an $A$-order of $k(\pi)$ containing $\pi$ and such that $\mathcal{O}$ is maximal at the unique zero $v_0$ of $\pi$ in $k(\pi)$ lying over the place $v$ of $k$.\\
If $~\mathcal{O}\otimes A_v \not\cong End\phi \otimes A_v $ then there exists a quotient Drinfeld module\\ $\psi = \phi / G_{\mathcal{L}}$ such that\\ $ End\psi \otimes A_v \cong \mathcal{O}\otimes A_v $ and $End\psi \otimes A_{\omega} \cong End\phi \otimes A_{\omega} $ at all the other places $\omega \neq v$ of $k$.
\end{lem}
\begin{pf}
With the hypothesis of the lemma,\\
we assume that $End\phi \otimes A_v \not\cong \mathcal{O}\otimes A_v$ as $A_v$-modules. That means there must exist at least one other place $v_1 \neq v_0$ of $k(\pi)$ lying over the place $v$ of $k$ (i.e. $\phi$ is not supersingular) such that the completion $\mathcal{O}_{v_1}$ of $\mathcal{O}$ at the place $v_1$ is different from the completion $\left(End\phi\right)_{v_1}$ of $End\phi$ at that same place $v_1$.\\
Let $v_0, ~v_1, \cdots, v_s$ be the places of $k(\pi)$ lying over the place $v$ of $k$. We choose $v_0$ here to be the unique zero of $\pi$ in $k(\pi)$ lying over the place $v$.\\ 
We are looking for a quotient Drinfeld module $\psi = \phi/G_{\mathcal{L}}$ such that\\ $End\psi\otimes A_v \cong \mathcal{O}\otimes A_v$ and $End\psi \otimes A_{\omega} \cong End\phi \otimes A_{\omega}$ at all the other places $\omega \neq v$.\\
The idea here is to act on the \'etale part of the Dieudonn\'e module $T_v\phi$ of $\phi$ so that the resulting endomorphism ring meets our needs.\\
Let then $M(x)$ be the minimal polynomial (Weil polynomial) of $\pi$ over $k$. \\
We know that the places $v_0, ~v_1, \cdots, v_s$ are described by the irreducible factors of $M(x)$ in $k_v[x]$. Let then $M(x)=M_0(x) \cdot M_1(x) \cdots M_s(x) \in k_v[x]$ be the irreducible decomposition of $M(x)$ over the completion field $k_v$.\\
We also know that the irreducible factor $M_0(x)=:M_{loc}(x)$ describing the zero $v_0$ of $\pi$ in $k(\pi)$ is the characteristic polynomial of the action of $\pi$ on the local part of the Dieudonn\'e module $\left(T_v\phi\right)_{loc}$ (see corollary~\ref{Mloc}).\\ In addition, $M_0(x) \equiv x^h \mod \mathfrak{p}_v$, where $h$ is the height of $\phi$ (see lemma~\ref{heightDM}). $M_{\text{\'et}}(x)=M_1(x) \cdots M_s(x)$ is the characteristic polynomial of the action of $\pi$ on the \'etale part of the Dieudonn\'e module $\left(T_v\phi\right)_{\text{\'et}}$. In this case, we therefore clearly see that $rank_W\left(T_v\phi\right)_{\text{\'et}}=\deg M_{\text{\'et}}(x) \geq 2$. Because if we had $\deg M_{\text{\'et}}(x) = 0$, $\phi$ would be supersingular and if we had $\deg M_{\text{\'et}}(x)=1$, $End\phi\otimes A_v$ and $\mathcal{O}\otimes A_v$ would be both maximal orders of the $k_v$-algebra $k_v(\pi)$ and thus we would have\\ $End\phi\otimes A_v \cong \mathcal{O}\otimes A_v$, which in either case contradicts our assumption. \\
We recall the notation $K_v$ which is the unique degree $m$ unramified extension of $k_v$ and $W$ its ring of integers. \\
We know that $\mathcal{O} \otimes A_v = \displaystyle\prod_{v_i \mid v}\mathcal{O}_{v_i}$ is an $A_v$-order of the $k_v$-algebra\\ $k_v(\pi)=End\phi\otimes k_v \cong End_{K_v[F,V]}V_v\phi$ (see remark~\ref{rmkWaterhouse}).\\
i.e. $\mathcal{O}\otimes A_v \subseteq k_v(\pi) \cong End_{K_v[F,V]}V_v\phi$. \\
Also, $\mathcal{O}$ is maximal at $v_0$ i.e. the completion $\mathcal{O}_{v_0}$ is the maximal order of the field $k_v(\pi_0)=k_v[x]/M_0(x)\cdot k_v[x]$.\\
Thus there exists a $W$-lattice $\mathcal{L}_0$ of $\left(V_v\phi\right)_{\text{\'et}}=\left(T_v\phi
\right)_{\text{\'et}} \otimes K_v$ containing $\left(T_v\phi\right)_{\text{\'et}}$ stable under the actions of $F$ and $V$, \\
(i.e. $T_v\phi = \left(T_v\phi\right)_{loc} \oplus \left(T_v\phi\right)_{\text{\'et}} \subseteq \left(T_v\phi \right)_{loc} \oplus \mathcal{L}_0 \subseteq V_v\phi=\left(V_v\phi\right)_{loc} \oplus \left(V_v\phi\right)_{\text{\'et}}$)\\
 such that the corresponding order $End_W\left(\left(T_v\phi\right)_{loc} \oplus \mathcal{L}_0 \right) \cong \mathcal{O}\otimes A_v$.\\
 We set $l=r-h=\deg M_{\text{\'et}}(x) \geq 2$. Let $(t_1, \cdots, t_l)$ be a $W$-basis of $\left(T_v\phi \right)_{\text{\'et}}$ and $(z_1, \cdots, z_l)$ be a $W$-basis of $\mathcal{L}_0$. $N_0$ denotes the matrix in $\mathscr{M}_{l \times l}\left(W \right)$ such that 
 $$
 \begin{pmatrix}
 t_1 \\ \vdots \\ t_l
 \end{pmatrix}= N_0
 \begin{pmatrix}
 z_1 \\ \vdots \\ z_l
 \end{pmatrix}
 $$  
 Let $s_0=v\left(det N_0 \right)$. Since $End\phi \otimes A_v \cong End_{W[F,V]}T_v\phi \not\cong \mathcal{O} \otimes A_v$, $s_0 \geq 1$.\\
 Since $K_v$ is an unramified extension of $k_v$ and the corresponding ring of integers $W$ is a discrete valuation ring, we keep (by abuse of language) the same notation $v$ for the place of $K_v$ extending the place $v$ of $k_v$. $\mathfrak{p}_v$ denotes the corresponding prime.\\
 $\det N_0= \beta_0 \mathfrak{p}_v^{s_0}$, where $\beta_0$ is a unit in $W$. The same way we did before, let us consider the morphism
 $$
 \begin{tabular}{rcl}
 $Co(N_0)^t: \left(T_v\phi\right)_{\text{\'et}}$ & $ \xlongrightarrow{\hspace{2cm}} $ & $\mathcal{L}_0$\\
 $\begin{pmatrix}
 t_1 \\ \vdots \\ t_l
 \end{pmatrix}$ & $\longmapsto$ & 
 $Co(N_0)^t
 \begin{pmatrix}
 t_1 \\ \vdots \\ t_l
 \end{pmatrix}=\beta_0\mathfrak{p}_v^{s_0}
\begin{pmatrix}
 z_1 \\ \vdots \\ z_l
 \end{pmatrix} 
 $
 \end{tabular}
 $$
where $Co(N_0)^t$ denotes the transpose of the co-matrix of $N_0$. We recall that $Co(N_0)^t \cdot N_0 = \det N_0 \cdot IdentityMatrix$.\\
The kernel of $Co(N_0)^t$ is given by $Ker\left(Co(N_0)^t\right) = N_0 \cdot D\left( \phi[\mathfrak{p}_v^{s_0}]\right)_{\text{\'et}}$. \\
$D\left( \phi[\mathfrak{p}_v^{s_0}]\right)_{\text{\'et}}$ is the $W[F,V]$-module associated to the group-scheme $\phi[\mathfrak{p}_v^{s_0}]_{\text{\'et}}$ (see remark~\ref{Dieudonne-rmk}). Indeed,\\
Let $\lambda_1t_1 + \cdots + \lambda_lt_l \in N_0 \cdot D\left( \phi[\mathfrak{p}_v^{s_0}]\right)_{\text{\'et}}$. We have then, \\
$Co(N_0)^t\cdot \left(\lambda_1t_1 + \cdots + \lambda_lt_l\right) \in Co(N_0)^t \cdot N_0 \cdot D\left( \phi[\mathfrak{p}_v^{s_0}]\right)_{\text{\'et}}= \mathfrak{p}_v^{s_0}\cdot D\left( \phi[\mathfrak{p}_v^{s_0}]\right)_{\text{\'et}}$=$\{0\}$.\\
We recall that $D\left(\phi[\mathfrak{p}_v^n]\right)$ can be identified to $T_v\phi / \mathfrak{p}_v^n\cdot T_v\phi$ for any $n\in \mathbb{N}$.\\
Conversely, let $\lambda_1t_1 + \cdots + \lambda_lt_l \in Ker\left(Co(N_0)^t\right)$ i.e. $Co(N_0)^t\cdot \left(\lambda_1t_1 + \cdots + \lambda_lt_l\right)$=$0$\\
That means $\beta_0\mathfrak{p}_v^{s_0}(\lambda_1z_1 + \cdots + \lambda_lz_l)=0$ and then\\ 
$\lambda_1z_1 + \cdots + \lambda_lz_l \in D\left(\phi[\mathfrak{p}_v^{s_0}]\right)_{\text{\'et}}$.\\ 
But $\lambda_1t_1 + \cdots + \lambda_lt_l=N_0\cdot \left(\lambda_1z_1 + \cdots + \lambda_lz_l\right) \in N_0\cdot D\left(\phi[\mathfrak{p}_v^{s_0}]\right)_{\text{\'et}}$.\\
Therefore $Ker\left(Co(N_0)^t\right)= N_0\cdot D\left(\phi[\mathfrak{p}_v^{s_0}]\right)_{\text{\'et}}$.\\
Applying the first isomorphism theorem to our morphism, one gets that \\
$\left(T_v\phi\right)_{\text{\'et}}/N_0\cdot D\left(\phi[\mathfrak{p}_v^{s_0}]\right)_{\text{\'et}} \cong Im\left(Co(N_0)^t\right)=\langle \mathfrak{p}_v^{s_0}z_1, \cdots , \mathfrak{p}_v^{s_0}z_l\rangle$.\\
Let $\mathcal{L}_{s_0}$ be the $W$-lattice generated by $\left(\mathfrak{p}_v^{s_0}z_1, \cdots , \mathfrak{p}_v^{s_0}z_l \right)$.\\ i.e. $\left(T_v\phi\right)_{\text{\'et}}/N_0\cdot D\left(\phi[\mathfrak{p}_v^{s_0}]\right)_{\text{\'et}} \cong \mathcal{L}_{s_0}$.\\
$N_0\cdot D\left(\phi[\mathfrak{p}_v^{s_0}]\right)_{\text{\'et}}$ is stable under the actions of $F$ and $V$ because $N_0$ commutes with the actions of $F$ and $V$ (via the stability of $\left(T_v\phi\right)_{\text{\'et}}$ and $\mathcal{L}_0$ under those actions) and $D\left(\phi[\mathfrak{p}_v^{s_0}]\right)$ is by definition stable under those actions (see theorem~\ref{DieudonneCartier}). \\
Let $G_{s_0}$ be the finite commutative $L$-group scheme associated to the $W[F,V]$-module $N_0\cdot D\left(\phi[\mathfrak{p}_v^{s_0}]\right)_{\text{\'et}}$ (theorem~\ref{DieudonneCartier}). We consider the additive separable polynomial 
$$u=x\displaystyle\prod_{\substack{\alpha \in G_{s_0} \\ \alpha \neq 0}}\left(1-\frac{x}{\alpha}\right)$$
whose kernel is $G_{s_0}$. By definition, $G_{s_0}$ is stable under the action of $\pi=F^m$. For the same reason as the case $\omega \neq v$ in lemma~\ref{lemma1}, $u \in L\{\tau\}$ and $u$ is an isogeny from the Drinfeld module $\phi$ to a Drinfeld module $\psi$. That is, $\phi_T \cdot u = u \cdot \psi_T$. In fact $\psi:= \phi /G_{s_0}$.\\
The Dieudonn\'e module of $\psi$ is given as follows:\\
$T_v\psi=T_v\left(\phi/G_{s_0}\right) \cong T_v\phi / D(G_{s_0})=\left(\left(T_v\phi\right)_{loc} \oplus \left(T_v\phi\right)_{\text{\'et}}\right)/N_0\cdot D\left(\phi[\mathfrak{p}_v^{s_0}]\right)_{\text{\'et}}$.\\
That is,\\
$T_v\psi \cong \left(T_v\phi\right)_{loc} \oplus \left(T_v\phi\right)_{\text{\'et}}/N_0\cdot D\left(\phi[\mathfrak{p}_v^{s_0}]\right)_{\text{\'et}}
\cong \left(T_v\phi\right)_{loc} \oplus \mathcal{L}_{s_0} $   \\
One easily checks that since $\mathcal{L}_{s_0}= \mathfrak{p}_v^{s_0}\cdot \mathcal{L}_0,~~ End_W\mathcal{L}_{s_0}=End_W\mathcal{L}_0$.\\
Therefore $End_W\left(\left(T_v\phi\right)_{loc} \oplus \mathcal{L}_{s_0} \right) \cong End_W\left(\left(T_v\phi\right)_{loc} \oplus \mathcal{L}_0 \right) \cong \mathcal{O}\otimes A_v$ and from the stability under the actions of $F$ and $V$, one gets\\
$\mathcal{O}\otimes A_v \cong End_{W[F,V]}\left(\left(T_v\phi\right)_{loc} \oplus \mathcal{L}_0 \right) \cong End_{W[F,V]}\left(\left(T_v\phi\right)_{loc} \oplus \mathcal{L}_{s_0} \right) \cong End_{W[F,V]}T_v\psi$\\
Hence $\mathcal{O}\otimes A_v \cong End \psi \otimes A_v$ (Thanks to the Tate's theorem~\ref{Tatetheorem2}).\\
At all the other places $\omega \neq v$ we have the exact sequence
$$0 \longrightarrow G_{s_0} \xhookrightarrow{\hspace{1cm}} \phi \xlongrightarrow{\hspace{0.4cm} u \hspace{0.4cm}} \psi=\phi/G_{s_0} \longrightarrow 0$$
Applying the Tate's theorem at the place $\omega$, we get 
$$ 0 \longrightarrow T_{\omega}\phi \longrightarrow T_{\omega}\psi  \longrightarrow 0$$
In fact by definition of the Dieudonn\'e functor in theorem~\ref{DieudonneCartier} and from the Lagrange theorem for finite group scheme, we have the following:\\
If $r_0=rank\left( N_0 \cdot D\left(\phi[\mathfrak{p}_v^{s_0}]\right)\right)$ then $\mathfrak{p}_v^{r_0} \cdot G_{s_0}= \{0\}$ i.e. $G_{s_0} \subseteq \phi[\mathfrak{p}_v^{r_0}]$. That means the Tate module $T_{\omega}G_{s_0}=\{0\}$ for any place $\omega \neq v$.\\
Hence we get from the above exact sequence that $T_{\omega} \phi \cong T_{\omega}\psi$.\\ In other words\\
$End\phi \otimes A_{\omega} \cong End_{A_{\omega}[\pi]}T_{\omega}\phi \cong End_{A_{\omega}[\pi]}T_{\omega}\psi \cong End\psi \otimes A_{\omega}$. 
\hspace{13cm}
\end{pf}
~\\
\noindent Let us recall the theorem we want to prove.

\begin{thm}\label{endo_thm}
$A=\mathbb{F}_q[T],~ k=\mathbb{F}_q(T)$ and $\mathfrak{p}_v$ is the (generator of the) kernel of the characteristic morphism $\gamma: A \longrightarrow L$ defining the finite $A$-field $L$.\\
$M(x)= x^r + a_1x^{r-1} + \cdots + a_{r-1}x + \mu\mathfrak{p}_v^m \in A[x]$ is a Weil polynomial, where $m=[L: A/\mathfrak{p}_v\cdot A]$. Let $\mathcal{O}$ be an $A$-order of the function field\\ $k(\pi)=k[x]/M(x)\cdot k[x]$. Let $v_0$ be the unique zero of $\pi$ in $k(\pi)$ lying over the place $v$ of $k$.\\
$\mathcal{O}$ is the endomorphism ring of a Drinfeld module in the isogeny class defined by the Weil polynomial $M(x)$ if and only if $\mathcal{O}$ contains $\pi$ and $\mathcal{O}$ is maximal at the place $v_0$. 
\end{thm}

\begin{pf}
With the hypotheses of the theorem, we have the following:\\
\fbox{\parbox[b]{0.5cm}{$\Rightarrow$}} If $\mathcal{O}=End\phi$ then it is clear that $\mathcal{O}$ contains the Frobenius endomorphism $\pi$. Yu proved in~\cite{yu1995isogenies} that $End\phi$ is maximal at the  zero $v_0$ of $\pi$ in $k(\pi)$.\\
~\\
\fbox{\parbox[b]{0.5cm}{$\Leftarrow$}} Conversely, let us assume that $\mathcal{O}$ contains $\pi$ and $\mathcal{O}$ is maximal at the place $v_0$.\\
Let $\phi$ be any Drinfeld module over $L$ in the isogeny class defined by $M(x)$.\\
We know that $\mathcal{O}$ and $End\phi$ differ at only finitely many places, since both are orders of the same function field $k(\pi)$. That means there exist finitely many places $\omega_1, \cdots, \omega_s$ such that \\
$\mathcal{O} \otimes A_{\omega} \cong End\phi \otimes A_{\omega}$ for all places $\omega$ except (may be) at $\omega \in \{ v,~\omega_1,~\omega_2, \cdots , \omega_s\}$.\\
For $\omega=\omega_1$, one can get from lemma~\ref{lemma1} a Drinfeld module $\phi_1$ defined over $L$ such that\\ $End\phi_1 \otimes A_{\omega_1} \cong \mathcal{O} \otimes A_{\omega_1}$ and\\
 $End\phi_1 \otimes A_{\nu} \cong End \phi \otimes A_{\nu}$ at all other places $\nu \neq \omega_1, v$. \\
That means $End\phi_1 \otimes A_{\omega} \cong \mathcal{O} \otimes A_{\omega}$ for all places $\omega$ of $k$ except (may be) at $\omega \in \{ v,~ \omega_2,~ \omega_3, \cdots, \omega_s\}$.\\
Repeating the process at all the places $\omega_i$, one gets from lemma~\ref{lemma1} a Drinfeld module $\varphi$ defined over $L$ such that \\
$End\varphi \otimes A_{\omega} \cong \mathcal{O} \otimes A_{\omega}$ for all places $\omega$ of $k$ with $\omega \neq v$.\\
Concerning the place $v$, we know in addition that $\mathcal{O}$ is maximal at the unique zero $v_0$ of $\pi$ in $k(\pi)$ lying over the place $v$.\\
We can therefore apply lemma~\ref{lemma2} and get the following:
\begin{itemize}
\item If $\varphi$ (equivalently our isogeny class) is supersingular, then we already have $End\varphi \otimes A_v \cong \mathcal{O} \otimes A_v$ as maximal order of the $k_v$-algebra (which is actually in this case a field) $k_v(\pi)$.
\item If $\varphi$ (equivalently our isogeny class) is not supersingular and\\ $End\varphi \otimes A_v \not\cong \mathcal{O} \otimes A_v$, then there exists (see lemma~\ref{lemma2}) a Drinfeld module $\psi= \varphi / G_{\mathcal{L}}$ such that \\
$End\psi \otimes A_v \cong \mathcal{O} \otimes A_v$ and \\
$End\psi \otimes A_{\omega} \cong End\varphi \otimes A_{\omega} \cong \mathcal{O} \otimes A_{\omega}$ at all the other places $\omega \neq v$.\\
In any case, we get a Drinfeld module $\psi$ such that \\
$End \psi \otimes A_{\omega} \cong \mathcal{O} \otimes A_{\omega}$ at all the places $\omega$ of $k$.\\
Hence $\mathcal{O}= End\psi$. 
\end{itemize}   
\hspace{13cm}
\end{pf}

\section{Application: Endomorphism rings in some isogeny classes of rank 3 Drinfeld modules}
We give in this part a more specific description of the orders occurring as endomorphism of a Drinfeld module in the special case of an isogeny class of rank 3 Drinfeld modules.
As a direct consequence of theorem~\ref{endo_thm}, we have the following:
\begin{prop}\label{rank3-Endo} We keep the same notation we have in the above mentioned theorem.\\
Let $M(x)= x^3 + a_1x^2 + a_2x +\mu\mathfrak{p}_v^m$ be a rank 3 Weil polynomial.
\begin{enumerate}
\item[1)] If $\mathfrak{p}_v \nmid a_2$ then an $A$-order $\mathcal{O}$ of $k(\pi)$ is the endomorphism ring of a Drinfeld module in the isogeny class defined by $M(x)$ if and only if it contains the Frobenius $\pi \in \mathcal{O}$.
\item[2)] Otherwise (i.e. if $\mathfrak{p}_v \mid a_2$),
an order $\mathcal{O}$ of $k(\pi)$ occurs as endomorphism ring of a Drinfeld module in the isogeny class defined by $M(x)$ if and only if the Frobenius endomorphism $\pi \in \mathcal{O}$ and $\mathcal{O}$ is maximal at all the places of $k(\pi)$ lying over $v$ (i.e. $\mathcal{O} \otimes A_v$ is a maximal order of the $k_v$-algebra $k_v(\pi)$).
\end{enumerate}
\end{prop}
\begin{pf}
\begin{enumerate}
\item[1)] If $\mathfrak{p}_v \nmid a_2$ then $M(x) \equiv x(x^2+a_1x+a_2) \mod \mathfrak{p}_v$. That means (see corollary~\ref{Mloc}) the irreducible factor $M_{loc}(x)$ of $M(x)$ in $k_v[x]$ describing the unique zero $v_0$ of $\pi$ in $k(\pi)$ is a degree $1$ polynomial. Therefore any $A$-order of $k(\pi)$ containing $\pi$ is already maximal at $v_0$. The statement follows then from theorem~\ref{endo_thm}.
\item[2)] If $\mathfrak{p}_v \mid a_2$ then we have two sub-cases.
\begin{itemize}
\item If $\mathfrak{p}_v \nmid a_1$ then $M(x) \equiv x^2(x+a_1) \mod \mathfrak{p}_v$.\\
That means there are two places of $k(\pi)$ lying over the place $v$. The zero $v_0$ of $\pi$ which is described by the irreducible factor $M_{loc}(x)$ of $M(x)$ in $k_v[x]$ fulfilling $M_{loc}(x) \equiv x^2 \mod \mathfrak{p}_v$ (see corollary~\ref{Mloc}), and another place $v_1$ described by the irreducible factor $M_1(x)$ of $M(x)$ in $k_v[x]$ fulfilling $M_1(x) \equiv x+a_1 \mod \mathfrak{p}_v$.\\
As a consequence, $\deg M_1(x) = 1$. That means the completion of any $A$-order $\mathcal{O}$ of $k(\pi)$ containing $\pi$ at the place $v_1$ must be maximal. \\
It follows that, $\mathcal{O}$ is maximal at the zero $v_0$ of $\pi$ if and only if $\mathcal{O}$ is maximal at all the places ($v_0$ and $v_1$) of $k(\pi)$ lying over $v$ and the statement follows.
\item If $\mathfrak{p}_v \mid a_1$ then $M(x) \equiv x^3 \mod \mathfrak{p}_v$. That means the isogeny class defined by $M(x)$ is supersingular. In other words there is a unique place (the zero $v_0$ of $\pi$) of $k(\pi)$ lying over $v$ and the statement follows from theorem~\ref{endo_thm}. 
\end{itemize}
\end{enumerate}
\hspace{13.5cm}
\end{pf}

\begin{rmk}[Recall]~\\
To check that $\mathcal{O}\otimes A_v$ is a maximal $A_v$-order in the $k_v$-algebra $k_v(\pi)$ one can just check that the norm of the conductor $\mathfrak{c}$ of $\mathcal{O}$ is not divisible by $\mathfrak{p}_v$. We recall that the norm of the conductor can be gotten from the relationship between the discriminant of the order $\mathcal{O}$ and the discriminant of the field $k(\pi)$.
$$ disc\left(\mathcal{O}\right) = N_{k(\pi)/k}\left(\mathfrak{c}\right)\cdot disc\left(k(\pi)\right)$$ 
\end{rmk}

In the upcoming part, we want to explicitly compute the maximal order of the cubic function field $k(\pi)$ and all the sub-orders occurring as endomorphism ring of a rank-3 Drinfeld module. 

\begin{prop}\cite[Corollary 5.2]{landquist2010explicit}~\\
Let $M_0(x)=x^3+c_1x +c_2$ be the standard form of the polynomial $M(x)=x^3+ a_1x^2 +a_2x + \mu Q$. Where $c_1$ and $c_2$ are like computed in the algorithm~\ref{newalgo}. Let $disc\left( M_0(x) \right)= \lambda\displaystyle\prod_{i=1}^l D_i^i$ be the square-free factorization of $disc\left( M_0(x) \right)$. \\
The discriminant of the function field $k(\pi)$ is given by $$disc\left( k(\pi) \right)=\lambda D \gcd(D_2D_4, c_2)^2 \text{ where } D=\displaystyle\prod_{i~odd} D_i,~~\lambda \in \mathbb{F}_q^*.$$
\end{prop} 
We will not give the proof in details since it has already been done in \cite{landquist2010explicit}. We just remind that the proof strongly relies on the fact that\\ $M_0(x)=x^3 + c_1x +c_2$ is given in the standard form. That is, for any prime element $\mathfrak{p} \in A$, $v_{\mathfrak{p}}\left( c_1 \right) < 2 \text{ or }  v_{\mathfrak{p}}\left( c_2 \right) < 3 $. This condition forces the valuation of the discriminant $v_{\mathfrak{p}}\left( disc\left(M_0(x) \right)\right)=v_{\mathfrak{p}}\left(-4c_1^3 - 27c_2^2 \right)$ to be bounded and leads to the following lemma.
\begin{lem}\cite[theorem 2]{llorente1983effective}~\\
Let $k(\pi)/k$ be a cubic function field defined by the irreducible polynomial $M_0(x)=x^3 + c_1x +c_2$ given in the standard form. Let $D_0= disc\left( M_0(x) \right)$ and $\Delta_0= disc\left( k(\pi) \right)$. For any prime $\mathfrak{p}$ of $k$ we have the the following:
\begin{itemize}
\item[(1)] $v_{\mathfrak{p}}\left( \Delta_0 \right)=2$ if and only if $v_{\mathfrak{p}}\left( c_1 \right) \geq v_{\mathfrak{p}}\left( c_2 \right) \geq 1$.
\item[(2)] $v_{\mathfrak{p}}\left( \Delta_0 \right)=1$ if and only if $v_{\mathfrak{p}}\left( D_0 \right)$ is odd.
\item[(3)] $v_{\mathfrak{p}}\left( \Delta_0 \right)=0$ otherwise. 
\end{itemize} 
\end{lem}
\begin{rmk}
The index of $\tilde{\pi}$ can therefore be computed using the fact that\\ $disc\left(M_0(x) \right) = ind(\tilde{\pi})^2 disc\left( k(\pi) \right)$ i.e. $$I:=ind(\tilde{\pi})=\sqrt{\frac{disc\left(M_0(x) \right)}{disc\left( k(\pi) \right)}}$$ We recall that $\tilde{\pi}$ and $\pi$ define the same function field $k(\pi)=k(\tilde{\pi})$.
\end{rmk}

\begin{prop}\cite[theorem 6.4]{landquist2010explicit}\label{pro1} and \cite[lemma 3.1]{scheidler2004algorithmic}~\\
Let $M_0(x)=x^3+ c_1x +c_2$ be the standard form of the Weil polynomial $M(x)=x^3+ a_1x^2 + a_2x + \mu Q.$ $\pi$ denotes a root of $M(x)$ and $\tilde{\pi}=\displaystyle\frac{\pi + \frac{a_1}{3}}{\gcd(g_1, g_2)}$ is a root of $M_0(x)$. Let $\displaystyle\omega_1= \alpha_1 + \tilde{\pi}$ and $\displaystyle\omega_2=\frac{\alpha_2 + \beta_2\tilde{\pi} + \tilde{\pi}^2}{I}$, where  $\alpha_1,~ \alpha_2$ and $\beta_2$ are elements of $A$. \\
$(1, \omega_1, \omega_2 )$ is an integral basis of the cubic function field $k(\pi)=k(\tilde{\pi})$ if and only if  
$
\begin{cases}
3\beta_2^2+c_1 \equiv 0 \mod I\\
\beta_2^3+c_1\beta_2+c_2 \equiv 0 \mod I^2\\
\alpha_2 \equiv -2\beta_2^2 \equiv 2c_1/3 \mod I
\end{cases}
$ 

\end{prop}

~\newline
\noindent \begin{pf}
The proof mainly relies on the following two facts:\begin{itemize}
\item $disc(1, \tilde{\pi}, \tilde{\pi}^2)=I^2disc\left( k(\pi)/k\right)$
\item For $\omega_2=\frac{\alpha_2 + \beta_2\tilde{\pi} + \tilde{\pi}^2}{I}$ to be integral it is necessary that $$\omega_2^2=\frac{\left(\alpha_2 + \beta_2\tilde{\pi} + \tilde{\pi}^2\right)^2}{I^2} \text{ and } (\alpha_1 + \tilde{\pi})\omega_2 \text{ both lie in } A[1, \tilde{\pi}, \omega_2]$$ In other words there exist $\lambda_0, \mu_0,~ \lambda_1, \mu_1$ and $\lambda_2, \mu_2 \in A$ such that $\omega_2^2= \lambda_0 + \lambda_1\tilde{\pi} + \lambda_2\omega_2$ and $\tilde{\pi}\omega_2= \mu_0 + \mu_1\tilde{\pi} + \mu_2\omega_2$.
\end{itemize}
\hspace{14cm} 
\end{pf} 
\begin{cor}
$\alpha_1$ in the previous proposition can be assumed to be 0 because if $\displaystyle\left(1, \alpha_1 + \tilde{\pi}, \frac{\alpha_2 + \beta_2\tilde{\pi} + \tilde{\pi}^2}{I} \right)$ is an integral basis, then so is $\displaystyle\left(1, \tilde{\pi}, \frac{\alpha_2 + \beta_2\tilde{\pi} + \tilde{\pi}^2}{I} \right)$.
 \end{cor}
\noindent This is simply due to the fact that both triples have the same discriminant.
\begin{rmk}
One can therefore, given an isogeny class of Drinfeld modules described by the Weil polynomial $M(x)=x^3+ a_1x^2 + a_2x + \mu Q$, compute the corresponding maximal order $\mathcal{O}_{max}$ which is the $A$-module generated by $(1,~ \omega_1,~ \omega_2 )$ as mentioned before.
\end{rmk} 
Let $\mathcal{O}_{max}=\langle 1,~ \omega_1,\omega_2 \rangle=\Big\{ (X, Y, Z) 
\begin{pmatrix}
1\\ \omega_1 \\ \omega_2
\end{pmatrix} \Big\vert~ X,~ Y,~ Z \in A \Big\}$.\\
We want now to give a complete list of sub-orders of $\mathcal{O}_{max}$ occurring as endomorphism rings of Drinfeld modules. We know from proposition~\ref{rank3-Endo} that this is equivalent to looking for sub-orders containing $\pi$ and whose conductor's norm (in case $\mathfrak{p}_v \mid a_2$) is relatively prime to $\mathfrak{p}_v$.\\
\noindent Let then $\mathcal{O}=\langle \tilde{\omega_0},~ \tilde{\omega_1},~ \tilde{\omega_2} \rangle $ be a sub-order of $\mathcal{O}_{max}.$\\
$1 \in \mathcal{O}.~$ That means one can write without loss of generality\\ $\mathcal{O}= \langle 1,~ \tilde{\omega_1},~ \tilde{\omega_2} \rangle =\Big\{ (\tilde{X}, \tilde{Y}, \tilde{Z})
\begin{pmatrix}
1\\ \tilde{\omega_1} \\ \tilde{\omega_2}
\end{pmatrix} \Big\vert \tilde{X},~\tilde{Y},~\tilde{Z} \in A \Big\}$\\
But $\tilde{\omega_1} \text{ and } \tilde{\omega_2} \in \mathcal{O}_{max}$. That means $$\tilde{\omega_1}=\tilde{\alpha_1} + \tilde{\beta_1}\omega_1 + \tilde{\gamma_1}\omega_2  \text{ and } \tilde{\omega_2}=\tilde{\alpha_2} + \tilde{\beta_2}\omega_1 + \tilde{\gamma_2}\omega_2$$ for some $\tilde{\alpha_i},~\tilde{\beta_i},~\tilde{\gamma_i} \in A~~i=1,2.$ In other words,\\
$
\begin{pmatrix}
1\\ \tilde{\omega_1} \\ \tilde{\omega_2}
\end{pmatrix}=
\begin{pmatrix}
1 & 0 & 0 \\
\tilde{\alpha_1} & \tilde{\beta_1} & \tilde{\gamma_1}\\
\tilde{\alpha_2} & \tilde{\beta_2} & \tilde{\gamma_2}
\end{pmatrix}
\begin{pmatrix}
1\\ \omega_1 \\ \omega_2
\end{pmatrix}
$. Let 
$ M=
\begin{pmatrix}
1 & 0 & 0 \\
\tilde{\alpha_1} & \tilde{\beta_1} & \tilde{\gamma_1}\\
\tilde{\alpha_2} & \tilde{\beta_2} & \tilde{\gamma_2}
\end{pmatrix} \in \mathscr{M}_3(A)
$\\
Where $\mathscr{M}_3(A)$ denotes the ring of $3\times 3$ -matrices with entries in A.\\
$M$ can be transformed into the so-called Hermite normal form. That means there the exists a matrix $U \in GL_3(A)$ and an upper triangular matrix $H$ such that $U\cdot M=H.$\\
Some simple row operations show that the Hermite normal form of $M$ looks like 
\begin{equation}\label{equat1}
H=\begin{pmatrix}
1 & 0 & 0 \\
0 & c & b\\
0 & 0 & a
\end{pmatrix} \text{ with } \deg_T(b) < \deg_T(a)
\end{equation}
We therefore redefine $\tilde{\omega_1}$ and $\tilde{\omega_2}$ as $\tilde{\omega_1}=c\omega_1 + b \omega_2$ and $\tilde{\omega_2}=a\omega_2$.\\
The sub-lattice $\mathcal{O}$ can then be written as 
$$ \mathcal{O}=\langle 1,~\tilde{\omega_1},~\tilde{\omega_2} \rangle = \Big\{ (X,~Y,~Z)
\begin{pmatrix}
1 & 0 & 0 \\
0 & c & b\\
0 & 0 & a
\end{pmatrix}
\begin{pmatrix}
1\\ \omega_1\\ \omega_2
\end{pmatrix} \Big\vert X,~ Y,~Z \in A \Big\}$$ 
\begin{rmk}
One clearly notices that the sub-lattice $\mathcal{O}$ above is an order if and only if $\tilde{\omega_1}^2,~\tilde{\omega_2}^2 \text{ and } \tilde{\omega_1}\tilde{\omega_2}$ belong to $\mathcal{O}$
\end{rmk}
\noindent But  
$\begin{cases}
\tilde{\omega_1}^2=(c\omega_1 + b\omega_2)^2=c^2\omega_1^2 + 2bc\omega_1\omega_2 + b^2\omega_2^2\\
\tilde{\omega_2}^2=(a\omega_2)^2=a^2\omega_2^2\\
\tilde{\omega_1}\tilde{\omega_1}=(c\omega_1 + b\omega_2)(a\omega_2)=ac\omega_1\omega_2 + ab\omega_2^2
\end{cases}  
$\\
Thus 
$
\begin{pmatrix}
\tilde{\omega_1}^2\\ \tilde{\omega_2}^2 \\ \tilde{\omega_1}\tilde{\omega_2}
\end{pmatrix}= \underbrace{
\begin{pmatrix}
c^2 & b^2 & 2bc\\
0 & a^2 & 0\\
0 & ab & ac
\end{pmatrix}}_{M_1}
\begin{pmatrix}
\omega_1^2\\ \omega_2^2 \\ \omega_1\omega_2
\end{pmatrix}
$\\
As we have seen in proposition~\ref{pro1} and its corollary, \\
$~~\omega_1=\tilde{\pi}$ and $\displaystyle \omega_2=\frac{\alpha_2 + \beta_2\tilde{\pi} + \tilde{\pi}^2}{I}$ where $\tilde{\pi}^3+c_1\tilde{\pi} +c_2=0$\\
One can therefore compute $\omega_1^2, ~ \omega_2^2 \text{ and } \omega_1\omega_2$ in terms of $\omega_1$ and $\omega_2$. One gets \\
$
\begin{pmatrix}
\omega_1^2 \\ \omega_2^2 \\ \omega_1\omega_2
\end{pmatrix}= \underbrace{
\begin{pmatrix}
X_{11} & X_{12} & X_{13}\\
X_{21} & X_{22} & X_{23}\\
X_{31} & X_{32} & X_{33} 
\end{pmatrix}}_{M_2}
\begin{pmatrix}
1\\ \omega_1\\ \omega_2
\end{pmatrix}
$ where\\
$\displaystyle X_{11}=-\alpha_2, ~ X_{12}=- \beta_2, ~X_{13}=I$\\
$\displaystyle X_{21}=\frac{\alpha_2\beta_2^2 - c_1\alpha_2 + 3\alpha_2^2 -2c_2\beta_2}{I^2}, ~X_{22}=\frac{-\beta_2^3 -c_1\beta_2 -c_2}{I^2}, ~X_{23}=\frac{\beta_2^2-c_1+2\alpha_2}{I}$\\
$\displaystyle X_{31}=\frac{\alpha_2\beta_2-c_2}{I}, ~ X_{32}=\frac{-\beta_2^2-c_1 +\alpha_2}{I} \text{ and } X_{33}=\beta_2.$\\
Therefore 
$$
\begin{pmatrix}
\tilde{\omega_1}^2\\ \tilde{\omega_2}^2 \\ \tilde{\omega_1}\tilde{\omega_2}
\end{pmatrix}=M_1M_2
\begin{pmatrix}
1\\ \omega_1 \\ \omega_2
\end{pmatrix}=M_1M_2H^{-1}
\begin{pmatrix}
1\\ \tilde{\omega_1}\\ \tilde{\omega_2}
\end{pmatrix}
$$
\begin{rmk}\label{rrmk}
$\mathcal{O}$ is an order if and only if $M_1M_2H^{-1} \in \mathscr{M}_3(A)$
\end{rmk}
Let us now investigate the orders occurring as endomorphism ring of a rank 3 Drinfeld module.\\
We know that in addition to the above mentioned condition, $\mathcal{O}=\langle 1,~ \tilde{\omega_1},~ \tilde{\omega_2}\rangle$ must contain the Frobenius $\pi$. In other words, there should exist $a_0,~b_0,~c_0 \in A$ such that\\ $\pi=a_0 + b_0 \tilde{\omega_1} + c_0 \tilde{\omega_2}$. But $$\omega_1= \tilde{\pi} \text{ and } \tilde{\pi}=\frac{\pi + \frac{a_1}{3}}{\gcd(g_1, g_2)}$$
Also $\tilde{\omega_1}=c\omega_1 + b\omega_2$ and $\tilde{\omega_2}=a\omega_2$. Therefore 
$$\displaystyle -\frac{a_1}{3} + \gcd(g_1, g_2)\cdot \omega_1= a_0 +b_0c\cdot \omega_1 +(b_0b + c_0a)\cdot \omega_2 $$ 
Thus 
\begin{equation}\label{equat2}
b_0c=\gcd(g_1, g_2) \text{ and } b_0b=-c_0a
\end{equation}
That is, 
$$
\begin{cases}
c \text{ divides } \gcd(g_1, g_2) \text{ and}\\
a \text{ divides } b\displaystyle\frac{\gcd(g_1, g_2)}{c}
\end{cases}
$$
We summarize our discussion in the following theorem:
\begin{thm}
$A=\mathbb{F}_q[T]$ and $k=\mathbb{F}_q(T)$\\
Let $M(x)=x^3 +a_1x^2 +a_2x + \mu Q \in A[x]$ be a Weil polynomial. In order to put $M(x)$ in a simple form $x^3 + b_1x + b_2$, let $b_1=\displaystyle\frac{-a_1^2}{3} + a_2 \text{ and }  b_2=\displaystyle\frac{2a_1^3}{27}-\frac{a_1a_2}{3} +\mu Q$ whose square-free factorizations are given by
$$ b_1=\mu_1\displaystyle\prod_{i=1}^{n_1} b_{1i}^i ~~ b_2= \mu_2\displaystyle\prod_{j=1}^{n_2} b_{2j}^j~~~ \mu_1, \mu_2 \in \mathbb{F}_q^* $$
In order to get the standard form $M_0(x)=x^3+c_1x+c_2$ of $M(x)$ (as defined in \ref{def_stdform}), we consider $g_1$ and $g_2$ the elements of $A$ defined by $$g_1=\displaystyle\prod_{i=1}^{n_1} b_{1i}^{\floor{\frac{i}{2}}} \text{ and } g_2= \displaystyle\prod_{j=1}^{n_2}b_{2j}^{\floor{\frac{j}{3}}}$$
We remove out from $b_1$ and $b_2$ resp. the highest square common divisor and the highest cubic common divisor by setting 
$$ c_1=\frac{b_1}{\gcd(g_1, g_2)^2} \text{ and } c_2=\frac{b_2}{\gcd(g_1, g_2)^3}$$
Let $\tilde{\pi}=\frac{\pi + \frac{a_1}{3}}{\gcd(g_1, g_2)}$ be a root of the standard polynomial $x^3+c_1x + c_2$.\\
Let $I=ind(\tilde{\pi} )= \frac{ind(\pi)}{gcd(g_1, g_2)^3}$, $\alpha_2 \text{ and } \beta_2 \in A $ such that\\
$$
\begin{cases}
3\beta_2^2+c_1 \equiv 0 \mod I\\
\beta_2^3+c_1\beta_2+c_2 \equiv 0 \mod I^2\\
\alpha_2 \equiv -2\beta_2^2 \equiv 2c_1/3 \mod I
\end{cases}
$$ 
We consider the matrix $M_2 \in \mathscr{M}_3\left( k \right)$ defined by \\
$$M_2=
\begin{pmatrix}
X_{11} & X_{12} & X_{13}\\
X_{21} & X_{22} & X_{23}\\
X_{31} & X_{32} & X_{33} 
\end{pmatrix}
\text{ where } $$
%
%
$\displaystyle X_{11}=-\alpha_2, ~ X_{12}=- \beta_2, ~X_{13}=I$\\
$\displaystyle X_{21}=\frac{\alpha_2\beta_2^2 - c_1\alpha_2 + 3\alpha_2^2 -2c_2\beta_2}{I^2}, ~X_{22}=\frac{-\beta_2^3 -c_1\beta_2 -c_2}{I^2}, ~X_{23}=\frac{\beta_2^2-c_1+2\alpha_2}{I}$\\
$\displaystyle X_{31}=\frac{\alpha_2\beta_2-c_2}{I}, ~ X_{32}=\frac{-\beta_2^2-c_1 +\alpha_2}{I} \text{ and } X_{33}=\beta_2.$\\

\noindent The Endomorphism rings of Drinfeld modules in the isogeny class defined by the Weil polynomial $M(x)$ are:
$$ \mathcal{O}= A + A\cdot \left( c\tilde{\pi} + b\left(\frac{\alpha_2 +\beta_2\tilde{\pi} + \tilde{\pi}^2}{I}\right) \right) + A\cdot a\left( \frac{\alpha_2 +\beta_2\tilde{\pi} + \tilde{\pi}^2}{I} \right)$$
such that $M_1M_2H^{-1} \in \mathscr{M}_3\left( A \right)$ and in addition $gcd(\mathfrak{p}_v,ac)=1$ if $\mathfrak{p}_v \mid a_2$. Where 
$$M_1=
\begin{pmatrix}
c^2 & b^2 & 2bc\\
0 & a^2 & 0\\
0 & ab & ac
\end{pmatrix} \text{ and }
H=\begin{pmatrix}
1 & 0 & 0 \\
0 & c & b\\
0 & 0 & a
\end{pmatrix}
$$
$
\begin{cases}
c \text{ runs through the divisors of } \gcd(g_1, g_2)\\
a \text{ runs through the divisors of } I\\
b \in A \text{ such that } \deg_T b < \deg_T a \text{ and } a \mid b\displaystyle\frac{\gcd(g_1, g_2)}{c}

\end{cases}
$
\end{thm}
\begin{pf}
The proof follows straightforwardly from our discussion before. The condition $gcd(\mathfrak{p}_v , ac)=1$ comes from the fact that in case $\mathfrak{p}_v \mid a_2$, the norm of the conductor of $\mathcal{O}$ must be prime to $\mathfrak{p}_v$ (see Proposition~\ref{rank3-Endo}). \hspace{7cm}
\end{pf}
\begin{cor}
Let $M(x)=x^3+a_1x^2 + a_2x + \mu Q \in A[x]$ be a Weil polynomial. $\pi$ is a root of $M(x)$ and $\tilde{\pi}=\displaystyle\pi + \frac{a_1}{3}$. Let \\
$b_1=\displaystyle\frac{-a_1^2}{3} + a_2$ and $b_2=\displaystyle \frac{2a_1^3}{27} - \frac{a_1a_2}{3} + \mu Q$.\\
If there is no prime $\mathfrak{p} \in A$ such that $\mathfrak{p}^2 \mid b_1$ and $\mathfrak{p}^3 \mid  b_2$ (in particular if $b_1$ and $b_2$ are coprime or $b_1$ is square-free or $b_2$ is cubic-free) then the endomorphism rings of Drinfeld modules in the isogeny class defined by the Weil polynomial $M(x)$ are 
$$\mathcal{O}_a = A + A\cdot\tilde{\pi} + A\cdot a\left(\frac{\alpha_2 + \beta_2\tilde{\pi} + \tilde{\pi}^2}{I}\right) $$
such that $M_aM_2H_a^{-1} \in \mathscr{M}_3(A)$ and in addition $gcd(\mathfrak{p}_v,a)=1$ if $\mathfrak{p}_v \mid a_2$. Where
$$ M_a=
\begin{pmatrix}
1 & 0 & 0 \\
0 & a^2 & 0 \\
0 & 0 & a
\end{pmatrix} \text{ and } H_a=
\begin{pmatrix}
1 & 0 & 0\\
0 & 1 & 0\\
0 & 0 & a
\end{pmatrix}
$$
Here $a$ runs through the divisors of the index $I=ind(\tilde{\pi})$.
\end{cor}
\begin{pf}
One can just reconsider the equation~(\ref{equat2}) right after remark~\ref{rrmk}. Here $\gcd( g_1, g_2)=1$. Thus $b_0c=1$ i.e. $b_0$ and $c$ are units. In addition $b_0b=-c_0a$ and $b_0$ is a unit. That means $a \mid b$. But $\deg_T b < \deg_T a$ (see equation~(\ref{equat1})). Therefore $b=0$. Hence the matrix $H$ in equation~(\ref{equat1}) and the matrix $M_1$ become
$$H_a=
\begin{pmatrix}
1 & 0 & 0\\
0 & 1 & 0\\
0 & 0 & a
\end{pmatrix} \text{ and } M_a=
\begin{pmatrix}
1 & 0 & 0\\
0 & a^2 & 0 \\
0 & 0 & a
\end{pmatrix}
$$ 
and the result follows. \hspace{9cm}
\end{pf}

\newpage
\bibliography{biblio}
\bibliographystyle{plain}
\end{document}